\documentclass[12pt]{amsart}
\usepackage{graphicx}       
\usepackage{amsmath}        
\usepackage{amssymb}        
\usepackage{amsfonts}       
\usepackage{amsthm}         
\usepackage[scr]{rsfso}     
\usepackage{setspace}	
\usepackage{comment}
\usepackage{mathtools}
\usepackage{enumerate}
\usepackage{helvet,courier}
\usepackage{mathptmx}
\usepackage{type1cm}
\DeclareMathAlphabet{\mathcal}{OMS}{cmsy}{m}{n}
\DeclareMathAlphabet{\mathbbold}{U}{bbold}{m}{n}  
\usepackage{multicol}
\usepackage{mathrsfs}
\usepackage{verbatim}
\usepackage[usenames,dvipsnames]{xcolor}

\theoremstyle{plain}
\newtheorem{theorem}{Theorem}
\newtheorem{lemma}[theorem]{Lemma}

\newcommand{\wh}{\widehat}

\let\S\relax

\newcommand{\R}{\mathbb{R}}
\newcommand{\N}{\mathbb{N}}
\newcommand{\Z}{\mathbb{Z}}

\newcommand{\S}{\mathscr{S}}

\DeclareMathOperator{\supp}{supp}

\let\l\relax
\let\r\relax

\newcommand{\l}{\left}
\newcommand{\r}{\right}
\newcommand{\abs}[1]{\l|#1\r|}
\newcommand{\norm}[1]{\l\|#1\r\|}
\DeclarePairedDelimiter\ceil{\lceil}{\rceil}

\title{Boundedness of bilinear radial Fourier multipliers}
\author{Petr Honz\'\i{}k and Maty\'a{}\v s Male\v cek}
\date{\today}

\thanks{ P. Honz\'\i{}k was supported by the grant GA\v{C}R P201/24-10505S.} 

\begin{document}

\begin{abstract}
We show that a bilinear  radial Fourier multiplier operator with symbol $\sigma$ is $L^2(\R^n)\times L^2(\R^n) \to L^1(\R^n)$ bounded, $n\in \mathbb N,$ if the function $\sigma$ satisfies the smoothness condition $\sigma(2^j\cdot)\Phi\in L^2_{1/2 +\epsilon}(\mathbb R^{2n})$ for  some $\epsilon>0$ and every $j\in \mathbb Z,$ where $\Phi$ is a smooth cutoff function adapted to the annulus $|x|\in [1/4,4]$. This condition is dimension free. We also apply similar reasoning to provide alternative proof of the initial result concerning multilinear Bochner-Riesz operator and prove an estimate for generalized bilinear Bochner-Riesz operator.       
\end{abstract} 

\maketitle

\section{Introduction}

A $m$-linear $(p_1,\dots, p_m,p)$ multiplier $\sigma(\xi_1,\dots,\xi_m)$ is a   function on $\mathbb R^{ n}\times \cdots \times \mathbb R^{ n}$
such that the corresponding $m$-linear operator  
$$
\begin{aligned}
&T_\sigma(f_1,\dots , f_m)(x)\\& = \int_{\mathbb R^{mn}}\sigma(\xi_1,\dots,\xi_m)\wh f_1(\xi_1)\cdots 
\wh f_m(\xi_m) e^{2\pi i x\cdot(\xi_1+\cdots+\xi_m)}d\xi_1\cdots d\xi_m, 
\end{aligned}
$$
 initially defined on $m$-tuples of Schwartz functions, 
has a  bounded extension 
from $L^{p_1}(\mathbb R^n)\times\cdots\times L^{p_m}(\mathbb R^n)$
to $L^{p}(\mathbb R^n)$ for appropriate $p_1,\dots, p_m,p$. 

It is known from the work in \cite{CM} for $p>1$ and  ~\cite{KS},~\cite{GT} for $p\le 1$, that the classical Mihlin
condition on $\sigma$ in $\mathbb R^{mn}$ yields boundedness for $T_\sigma$ from $L^{p_1}(\mathbb R^n)\times\cdots\times L^{p_m}(\mathbb R^n)$
to $L^{p}(\mathbb R^n)$ for all $1<p_1,\dots p_m\le \infty$, $1/m<p =  (1/p_1+\cdots +1/p_m)^{-1}<\infty$. 
The Mihlin condition in this setting is usually referred to as the Coifman-Meyer 
condition and  the associated multipliers bear the same names as well. 

In~\cite{GHH}, Grafakos, Honz\'\i{}k and He studied the sharpness of the estimate in the bilinear case. They 
focused on the  boundedness of $T_\sigma$ in case $L^2\times L^2\to L^1$  
and they proved the following sharp result.

\begin{theorem}\label{GHH} Suppose $\wh\psi\in\mathcal C_0^{\infty}(\mathbb R^{2n})$ is positive and  supported in the annulus 
$$\{(\xi,\eta):1/2\le|(\xi,\eta)|\le 2\}$$ such that $\sum_{j\in\mathbb Z}\wh\psi_j(\xi,\eta)=
\sum_{j}\wh\psi(2^{-j}(\xi,\eta))=1$ for all $(\xi,\eta)\neq0$.
Let $ 1< r<\infty$, $s>\max\{n/2, 2n/r\}$, and suppose  there is a constant $A$ such that
\begin{equation}\label{usbn}
\sup_{j}\|\sigma(2^j\cdot)\wh\psi\|_{L^r_s(\mathbb R^{2n})}\le A<\infty.
\end{equation}
Then there is a constant $C=C(n,\Psi)$ such that  the bilinear operator 
$$
T_{\sigma}(f,g)(x)=
\int_{\mathbb R^{2n}}\sigma(\xi,\eta)\wh f(\xi)\wh g(\eta)e^{2\pi i x\cdot (\xi+\eta)}d\xi d\eta,
$$
initially defined on   Schwartz functions $f$ and $g$,
satisfies   
\begin{equation}
\|T_{\sigma}(f,g)\|_{L^1(\mathbb R^n)}\le CA\|f\|_{L^2(\mathbb R^n)}\|g\|_{L^2(\mathbb R^n)}.
\end{equation}
\end{theorem}

Radial Fourier multipliers are an important class of operators, which includes operators such as the famous Bochner-Riesz multiplier. The study of radial Fourier multiplier operators has a long history and recently there many strong results were found, such as the famous result of Heo, Nazarov and Seeger~\cite{HNS}. 

The goal of this article is to initiate more focused study of general bilinear Fourier multipliers with radial symbols. We will focus on the case $L^2\times L^2\to L^1,$ where we are able to obtain good results via simple but novel method similar to that of Coifman and Meyer~\cite{CM}. We then also obtain a new simple proof for the boundedness of multilinear Bochner-Riesz in the case $L^2\times L^2\to L^1,$ and we explore further generalization of the result.      

First, we are going to show that when the symbol of the multiplier operator is radial, this result similar to~\cite{GHH} is in fact dimension free. Our main result is:

\begin{theorem}\label{Main} Suppose $\wh\psi\in\mathcal C_0^{\infty}(\mathbb R^{2n})$ is positive, radial and  supported in the annulus 
$$\{(\xi,\eta):1/4\le|(\xi,\eta)|\le 4\}$$ such that $\sum_{j\in\mathbb Z}\wh\psi_j(\xi,\eta)=
\sum_{j}\wh\psi(2^{-j}(\xi,\eta))=1$ for all $(\xi,\eta)\neq 0$. 
Suppose $\sigma$ is a radial function on $\mathbb R^{2n}$ such that there is $\epsilon>0$ and a constant $A$ such that
\begin{equation}\label{usb}
\sup_{j}\|\sigma(2^j\cdot)\wh\psi\|_{L^2_{1/2+\epsilon}(\mathbb R^{2n})}\le A<\infty.
\end{equation}
Then there is a constant $C$ such that  the bilinear operator 
$$
T_{\sigma}(f,g)(x)=
\int_{\mathbb R^{2n}}\sigma(\xi,\eta)\wh f(\xi)\wh g(\eta)e^{2\pi i x\cdot (\xi+\eta)}d\xi d\eta,
$$
initially defined on   Schwartz functions $f$ and $g$,
satisfies   
\begin{equation}
\|T_{\sigma}(f,g)\|_{L^1(\mathbb R^n)}\le CA\|f\|_{L^2(\mathbb R^n)}\|g\|_{L^2(\mathbb R^n)}.
\end{equation}
\end{theorem}  

As a by-product of the proof of this theorem, we also obtain two results concerning Bochner-Riesz operators. First, the following result on $m-linear$ Bochner-Riesz. This was proved in the bilinear case by Bernicot, Grafakos, Song and Yan~\cite{BGSY}, and in the multilinear case by He, Li and Zheng~\cite{HLZ}. Our proof is, however, simpler and self-contained. 

\begin{theorem} \label{Bochner_Riesz}
	Let $m \geq 2$ be a natural number and $\lambda > \frac{m}{2} - 1$. Then the $m$-linear Bochner-Riesz operator, given by multiplier
	$$
		\sigma^\lambda{(\xi_1, \dots, \xi_m)} =
		\l(
			1 - \l( \abs{\xi_1}^2 + \dots + \abs{\xi_m}^2 \r)
		\r)^\lambda_+, 
		\quad \xi_1, \dots, \xi_m \in \R^n
	$$
	is bounded from $\l(L^2(\R^n)\r)^m$ to $L^{\frac{2}{m}}(\R^n)$.
\end{theorem}

In the bilinear case, the bilinear Bochner-Riesz operator may be, using our method, easily generalized, we may for example prove:

\begin{theorem} \label{Modified_Bochner_Riesz}
	Let $\gamma > 1$. For $\lambda \geq 0$ define modified bilinear Bochner-Riesz operator as
	$$
		\sigma^{\lambda,\,\gamma}{(\xi_1, \xi_2)} =
		\frac{\sigma^\lambda(\xi_1, \xi_2)}{\l(1-\log{\l(1-\l[\abs{\xi_1}^2 + \abs{\xi_2}^2 \r]\r)}\r)^\gamma}
	$$
	for $ \abs{\xi_1}^2 + \abs{\xi_2}^2 < 1$ and 0 otherwise.

	Then $T_{\sigma^{0,\, \gamma}}$ is bounded from $L^2(\R^n) \times L^2(\R^n)$ to $L^1(\R^n)$.
\end{theorem}

This logarithmic generalization of Bochner-Riesz was considered in linear case by Seeger~\cite{Seeger_Modified_BR}. 

The method we used in this paper may be applied in many more situations,
such maximal or multilinear versions of the operators or other specific examples of radial multipliers. In some cases this leads to new results or new proofs of the best known results. We will address these in the second part of this paper, which is currently in preparation. 

We note that while for the case $L^2\times L^2 \to L^1$ the smoothness condition $s>1/2$ seems to be optimal, our method alone does not seem to give optimal results for other cases $L^{p_1}\times L^{p_2} \to L^p,$ where $p_1$ or $p_2$ is not equal to $2$. Therefore, in order to obtain optimal results for the full range new ideas will be needed.

\section{Definitions and lemmas}

By the Sobolev space $L^p_s\l(\R^n\r)$, where $p \geq 1$ and $s >0$, we understand the space of all tempered distributions $f \in \S'\l(\R^n\r)$ such that their norm
$\norm{( (1+|\cdot|^2)^\frac{s}{2} \hat{f} )^\vee }_p$ is finite. \\
For a fixed Schwartz function $\Psi \in \S\l(R^n\r)$, we define the Littlewood-Payley operator $\Delta^\Psi_j$ as
$\Delta^\Psi_j f = \l( \Psi_j\l(2^{-j}\cdot \r) \hat{f} \r)^\vee$.

Finally, if $\sigma$ is a function on $\R^n$, then we denote by $T_\sigma$ the operator given by the multiplier $\sigma$, i.e., $T_\sigma(f) = (\sigma \hat{f} )^\vee $ for $f \in \S(\R^n)$.

We call a function $f$ on $\mathbb R^{d}$ radial if it is either constant or undefined  on each sphere $t\mathcal S^{d-1},$  $0<t<\infty.$ We see that if a radial function $f$ is in any Sobolev space $L^p_{s},$ $s>0,$ then it is in the space $L^p$ and the trace function $\bar f(t)=f(te),$ where $e$ is an unit vector, is well defined and unique almost everywhere on $(0,\infty).$      

The following lemma is a version of the main Theorem~\ref{Main} formulated in terms of the trace. 

\begin{lemma}\label{main_lemma}
    Let $\sigma: \R^{2n} \to \R$ be a radial function and $\Phi \in \S{\l(\R^{2n}\r)}$ a radial function that is supported in an annulus $\frac{1}{4} \leq \abs{\xi} \leq 4$ and satisfies $\sum_{j \in \Z}{\Phi{\l(2^{-j} \xi\r)}} = 1$ for $\xi \neq 0$.
    
    Set $\sigma_j = \sigma \, \Phi(2^{-j}\,\cdot)$ and $\Tilde{\sigma}_j{\l(\pm \abs{\cdot}^2\r)} = \sigma_j{\l( 2^{j+3} \cdot \r)}$ for $j \in \Z$ and suppose, that for some $\varepsilon \in (0,\frac{1}{2})$, the following holds:
    $$
        C 
        \coloneqq
        \sup_{j \in \Z}{
            \norm{
                \Tilde{\sigma_j}
            }_{\frac{1}{2}+\varepsilon, \, 2}
        }
        < \infty.
    $$

    Then $\sigma$ is a bilinear multiplier from $L^2{\l(\R^n\r)} \times L^2{\l(\R^n\r)}$ to $L^1{\l(\R^n\r)}$
\end{lemma}

This gives the Theorem~\ref{Main} by the following lemma:

\begin{lemma}\label{trace}
Suppose a radial function $f$ is supported in an annulus $\frac 1 r<|x|<r,$ $r>1$ on $\mathbb R^d.$ Let $1\leq p<\infty,$ $s>1/p,$ $d\geq 2.$ 
Fix $x\in \mathbb R^d$ with $|x|=1.$ For $t>0$ denote $\bar f(t)=\tilde f(t^2) = f(tx).$ Suppose $f\in L^p_s(\mathbb R^d),$ then  $\bar f\in L^p_s(\mathbb R)$ and $\tilde f \in L^p_s(\mathbb R).$  
\end{lemma}

We get $\bar f\in L^p_s(\mathbb R)$ by the Theorem 5 from~\cite{SSV}. Then, as $t\mapsto t^2$ is a diffeomorphism  on the interval $[1/r,r]$ we also get $\tilde f\in L^p_s(\mathbb R).$

The Lemma~\ref{main_lemma} is sharp. It is well known example related to the Sobolev imbedding that there exists compactly supported function in $g\in L^2_{1/2}(\mathbb R)$ which is not in $L^\infty$. We may assume that the support is in the interval $[1,2]$ and then construct $f(x)=g(|x|^2)$ on $\mathbb R^{2n}.$ Then function $\sigma$ is not in $L^\infty,$ and therefore it is not a bilinear multiplier. 

To prove the Lemma~\ref{main_lemma}, we will need two additional lemmas. Lemma \ref{abs_convergence_lemma} gives us information about the behavior of a Fourier series of functions from $L^2_{\frac{1}{2}+\varepsilon}$.

\begin{lemma}\label{abs_convergence_lemma}
    Let $\varepsilon \in (0,\frac{1}{2})$ and $\sigma \in L^2_{\frac{1}{2}+\varepsilon}{\l(\R\r)}$ be a function supported in $[-\frac{1}{2}, \frac{1}{2}]$. Then $\sigma$ has an absolutely convergent Fourier series which uniformly converges to $\sigma$. 

	Furthermore, if for $j \in \N$ we denote $D_j = \l\{ k \in \Z: 2^{j-1} \leq \abs{k} < 2^j \r\}$ and $D_0 = \{0\}$, then there exists a $C>0$, independent of $\sigma$, such that the following inequality holds for any non-negative integer $j$:
    $$
        \sum_{k \in D_j}{
            \abs{\hat{\sigma}(k)}
        }
        \leq
        C \norm{\sigma}_{\frac{1}{2}+\varepsilon, \,2} \, 2^{-j\varepsilon}.
    $$
\end{lemma}

The key lemma is the Lemma \ref{key_lemma}, which gives us a decomposition of compactly supported radial multipliers. 

\begin{lemma}\label{key_lemma}
    Let $\sigma:\R^{mn} \to \R$ be a function supported in a ball $B_{\R^{mn}}(0,R)$ for some $R > 0$. Consider a function $\Tilde{\sigma}$ defined as
    $$
        \Tilde{\sigma}{\l( \pm \abs{\xi}^2\r)} 
        =
        \sigma{\l( \sqrt{2m} R \xi \r)}
    $$
	and suppose, that $\Tilde{\sigma}$ has absolutely convergent Fourier series, which converges to $\Tilde{\sigma}$.
    Then for any $f_1 \dots, f_m \in \S{\l( \R^n \r)}$ we have:
    $$
        T_\sigma{\l( f_1, \dots, f_m\r)}(x)
        =
        \sum_{k \in \Z}{
            c_k
            \prod_{j=1}^m{
                T_{\sigma^k}{\l( f_j\r)}(x)
            }
        },
    $$
    where $c_k$'s are Fourier coefficients of $\Tilde{\sigma}$ and  $\sigma^k$ are Schwartz functions, with supremum norm equal to 1, supported in a ball $B_{\R^n}(0,\sqrt{m}R)$. Specially, it is a Fourier multiplier on $L^2(\R^n)$ with operator norm less or equal to 1.
\end{lemma}

Lemma \ref{key_lemma} is essential in the proof of Lemma \ref{main_lemma}. We will also use it to prove previously known result, that is boundednes of m-linear Bochner-Riesz operator.

In both theorems \ref{Bochner_Riesz} and \ref{Modified_Bochner_Riesz}, we will use a slightly weaker version of theorem \ref{main_lemma}, which can be applied only for compactly supported multipliers. On the other hand, it yields boundedness not only for the bilinear case, but also for the $m$-linear case, where $m \geq 2$.

\begin{lemma}\label{compact_multipliers_lemma}
    Let $m \in \N$, $m \geq 2$ and $\sigma: \R^{mn} \to \R$ be a radial function supported in a ball $B_{\R^{mn}}(0,R)$ for some $R>0$. Let $\Tilde{\sigma}$ and $c_k$, where $k \in \Z$, be defined as in Lemma \ref{key_lemma}. Suppose that the sequence $\{c_k\}_{k\in\Z}$ lies in (possibly quasi-Banach) space $\ell^{\frac{2}{m}}(\Z)$. 
    
    Then $T_\sigma$ is bounded from $\l( L^2(\R^n) \r)^m$ to $L^{\frac{2}{m}}(\R^n)$ and for any $f_1, \dots,f_m \in L^2(\R^n)$ it holds, that
    $$
        \norm{
            T_\sigma(f_1,\dots,f_m)
        }_{2/m}
        \leq
        \norm{
            \{c_k\}_{k\in \Z}
        }_{\ell^{\frac{2}{m}}} \,
        \norm{f_1}_2
        \dots
        \norm{f_m}_2.
    $$
\end{lemma}

The proof of Lemma \ref{compact_multipliers_lemma} follows directly from Lemma \ref{key_lemma}, via H\"older inequality and summation.

\section{Proof of Lemma \ref{abs_convergence_lemma}}

By the results from \cite{H_Cont} we see that $L^2_{\frac{1}{2}+\varepsilon}(\R) \hookrightarrow C^{0, \varepsilon}(\R)$, so the Fourier series of $\sigma$ converges pointwise to $\sigma$. 

Now, let $\Psi_0 \in \S(\R)$ with support in a ball $\abs{\xi} \leq1$. Further set $\Psi \in \S(\R)$, which equals 1 in an annulus $\frac{1}{2} \leq \abs{\xi} \leq 1$ and $\Psi(0) = 0$. Then by \cite{Grafakos_Modern}:
\begin{equation}\label{eq_Frac_sobol+Triebel}
    \norm{\sigma}_{\frac{1}{2}+\varepsilon, \, 2}
    \approx
    \norm{
       \l(
            \abs{
                \l( 
                    \Psi_0 \hat{\sigma}
                \r)^\vee
            }^2
            +
            \sum_{j=1}^\infty{
                \l(
                    2^{j\l(\frac{1}{2}+\varepsilon\r)}
                    \abs{
                        \Delta_j^\Psi{\sigma}
                    }
                \r)^2
            }
       \r)^\frac{1}{2}
    }_2
\end{equation}

Take $j \in \N$. Then
$$
	\sum_{k \in D_j}{ \abs{\hat{\sigma}(k)}^2 }
	=
	\sum_{k \in D_j}{ \abs{\Psi\l(2^{-j}k\r)}^2 \abs{\hat{\sigma}(k)}^2 }
	\leq
	\norm{ \Delta_j^\Psi \sigma }_2^2.
$$
By using the Cauchy-Schwarz inequality and \eqref{eq_Frac_sobol+Triebel}, we get:
$$
	\sum_{k \in D_j}{ \abs{\hat{\sigma}(k)} }
	\leq
	2^{\frac{j}{2}} 
	\l(
		\sum_{k \in D_j}{ \abs{\hat{\sigma}(k)}^2 }
	\r)^\frac{1}{2}
	\leq
	2^{-\varepsilon j} \norm{ 2^{j\l( \frac{1}{2}+ \varepsilon \r) } \Delta_j^\Psi \sigma }_2
	\lesssim 
	2^{-\varepsilon j} 
	\norm{\sigma}_{\frac{1}{2}+\varepsilon, \, 2}.
$$
Now the fact that $L^2\left(\left[-\frac{1}{2},\frac{1}{2}\right]\right) \hookrightarrow L^1\left(\left[-\frac{1}{2}, \frac{1}{2}\right]\right)$ and that $L^2_{\frac{1}{2}+\varepsilon} {(\R)} \hookrightarrow L^2(\R)$, we see that there is a $C_1>0$ independent of $\sigma$, such that
$$
    \sum_{k \in D_0}{\abs{\hat{\sigma}(k)}}
    =
    \abs{\hat{\sigma}(0)}
    \leq
    \norm{\sigma}_1
    \leq
    C_1 \norm{\sigma}_{\frac{1}{2}+\varepsilon, \, 2}.
$$
From this, it follows that there is a $C > 0$ independent of $\sigma$, such that
$$
	\sum_{k \in \Z}{ \abs{\hat{\sigma}(k)} }
	=
	\sum_{j =0}^\infty{
		\l(
			\sum_{k \in D_j}{ \abs{\hat{\sigma}(k)} }
		\r)
	} \\
	\leq
	C 
    \l( 
        \sum_{j=0}^\infty{
            2^{-\varepsilon j}
        } 
    \r) 
    \norm{\sigma}_{\frac{1}{2}+\varepsilon, \, 2}.
$$
\qed

\section{Proof of key Lemma \ref{key_lemma}}

Notice that the support of $\Tilde{\sigma}$ lies in $\left[-\frac{1}{2}, \frac{1}{2}\right]$. We will again denote the Fourier coefficients of $\Tilde{\sigma}$ as $c_k$ for $k \in \Z$. We see that
    $$
        \sigma(\xi) 
        =
        \Tilde{\sigma}{
            \l(
                \frac{1}{2mR^2}
                \abs{\xi}^2
            \r)
        }
        =
        \sum_{k \in \Z}{
            c_k 
            e^{\frac{\pi}{m} k i \abs{\frac{1}{R}\xi}^2}
        }
    $$
    for $|\xi| \leq \sqrt{m} R$.

    Now choose a function $\phi \in \S{\l( \R^n \r)}$ such that $\phi(\xi) = 1$ for $|\xi| \leq 1$ and $\phi(\xi) = 0$ for $\abs{\xi} > \sqrt{m}$. Then, because $\sigma$ is supported in a $B_{\R^{mn}}(0,R)$, we have
    $$
        \sigma(\xi_1, \dots,\xi_m) 
        =
        \sum_{k \in \Z}{
            c_k 
            \prod_{j=1}^m{
                e^{\frac{\pi}{m} k i \abs{\frac{1}{R}\xi_j}^2}
                \phi{\l(\frac{1}{R}\xi_j\r)}
            }
        }        
        =
        \sum_{k \in \Z}{
            c_k 
            \prod_{j=1}^m{
                \sigma^k(\xi_j)
            }
        },    
    $$
    where we set 
    $\sigma^k(\xi_j) = e^{\frac{\pi}{m} k i \abs{\frac{1}{R}\xi_j}^2} \phi(\frac{1}{R}\xi_j)$. Using the absolute convergence of the sum of $c_k$'s concludes the proof of the lemma.
\qed

\section{Proof of Lemma \ref{main_lemma}}

    Pick $f_1, f_2 \in \S(\R^n)$ and choose some $\varphi \in \S(\R^n)$, which is supported in an annulus $\frac{1}{4} \leq |\xi| \leq 4$ such that $0 \leq \varphi\leq 1$ and $\sum_{j \in \Z}{\varphi(2^{-j}\xi)} = 1$ for $\xi \neq 0$.

    Observe that if $j, r, s \in \Z$ satisfy at least one of the following conditions: 
    \begin{itemize}
        \item $r > j+4$ or $s > j+4$,
        \item $r < j-4$ and $s < j-4$,
    \end{itemize}
    then $\sigma_j\l( \Delta^\varphi_r f_1 \otimes \Delta^\varphi_s f_2 \r)^\wedge = 0$:

    We use Lemma \ref{key_lemma} for each $j \in \Z$. Denote $c_{j,k}$ the Fourier coefficient of $\Tilde{\sigma}_j$ and $\sigma^k_j(\xi) = e^{\frac{\pi}{2} k i \abs{2^{-j-2} \xi}^2} \phi(2^{-j-2} \xi)$ for some $\phi \in \S(\R^n)$ that satisfies $\phi(\xi)=1$ for $\abs{\xi} \leq 1$, $\phi(\xi) = 0$ for $\abs{\xi} > \sqrt{2}$. Then
    $$
        T_{\sigma_j}(f_1,f_2) 
        =
        \sum_{k\in\Z}{
            c_{j,k} \, T_{\sigma_j^k}(f_1) \, T_{\sigma_j^k}(f_2)
        }.
    $$

    We can split the operator $T_\sigma$ into 3 parts: The diagonal part $T^{1,1}$ and two off-diagonal parts given by operator $T^{2}$, which are defined as:
    \begin{align*}
        T_\sigma^1 (f_1,f_2)
        &= 
        T_\sigma (f_1,f_2)
        - 
        \left(
                T_\sigma^2(f_1,f_2)
                +T_\sigma^2(f_2,f_1)
        \right),\\ 
        T_\sigma^2(f_1,f_2)
        &=
        \sum_{j\in\Z}{ \,
            \sum_{m =0}^\infty{ \,
                \sum_{k \in D_m}
                {
                    c_{j,k} \,
                    T_{\sigma_j^k}{
                    \l(
                        \sum_{r=-\infty}^{j-6-\ceil{\frac{m}{2}}}{\Delta_r^\varphi{f_1}}
                    \r)
                    }
                    T_{\sigma_j^k}{
                    \l(
                        \sum_{s=j-4}^{j+4}{\Delta_s^\varphi{f_2}}
                    \r)
                    }
                }
            }
        },
    \end{align*}
    where $D_m = \l\{ k \in \Z: \, 2^{m-1} \leq \abs{k} < 2^m\r\}$ for $m \in \N$ and $D_0 = \{0\}$. Notice that $T_\sigma(f_1,f_2)=T_\sigma^1(f_1,f_2)+ T_\sigma^2(f_1,f_2) + T_\sigma^2(f_2,f_1)$, so to prove boundedness of $T_\sigma$ it suffices to prove boundedness of $T_\sigma^1$ and $T_\sigma^2$.

    First, we will analyze the diagonal case and then we will focus on the off-diagonal parts.

 \subsection{Diagonal part} 

    Notice that because of the form of the support of the function $\sigma_j\l( \Delta^\varphi_r f_1 \otimes \Delta^\varphi_s f_2 \r)^\wedge$, we have
    \begin{multline*}
       T_\sigma^1
        = 
        \sum_{j\in\Z} \,
            \sum_{m =0}^\infty \,
                \sum_{k \in D_m}
                    c_{j,k} \,
                    \l(
                        T_{\sigma_j^k}{
                        \l(
                            \sum_{r=j-4}^{j+4}{\Delta_r^\varphi{f_1}}
                        \r)
                        }
                        T_{\sigma_j^k}{
                        \l(
                            \sum_{s=j-4}^{j+4}{\Delta_s^\varphi{f_2}}
                        \r)
                        } + \r. \\
                        \l.+
                        T_{\sigma_j^k}{
                        \l(
                            \sum_{r=j-5-\ceil{\frac{m}{2}}}^{j-5}{\Delta_r^\varphi{f_1}}
                        \r)
                        }
                        T_{\sigma_j^k}{
                        \l(
                            \sum_{s=j-4}^{j+4}{\Delta_s^\varphi{f_2}}
                        \r)
                        }
                        +
                        T_{\sigma_j^k}{
                        \l(
                            \sum_{r=j-4}^{j+4}{\Delta_r^\varphi{f_1}}
                        \r)
                        }
                        T_{\sigma_j^k}{
                        \l(
                            \sum_{s=j-5-\ceil{\frac{m}{2}}}^{j-5}{\Delta_s^\varphi{f_2}}
                        \r)
                        }
                    \r)
    \end{multline*}
    Consider an operator
    $$
        T^{1,1}_\sigma(f_1,f_2)
        =
        \sum_{m=0}^\infty{ \,
            \sum_{k \in D_m}{ \, \,
                \sum_{r,l =-5-\ceil{\frac{m}{2}}}^4{ \,
                    \sum_{j\in \Z}{ \,
                        \abs{c_{j,k}}
                        \abs{T_{\sigma_j^k}(\Delta_{j+r}^\varphi\,f_1)}
                        \abs{T_{\sigma_j^k}(\Delta_{j+l}^\varphi\,f_2)}
                    }
                }
            }
        }.
    $$
 The operator $T^{1,1}_\sigma$ clearly dominates $T^1_\sigma$, so it is sufficient to show boundedness only for $T_\sigma^{1,1}$.
    
    Using Fubini's theorem and repeatedly applying the Cauchy-Schwarz inequality for the sums and the integral, we get:
    $$
        \norm{
            T_\sigma^{1,1}(f_1,f_2)
        }_1
        \leq
        \sum_{m=0}^\infty{ \,
            \sum_{k\in D_m}{ \,
                \prod_{i=1}^2{ 
                    \l(
                        \sum_{l=-5-\ceil{m/2}}^4{ 
                            \l(
                                \sum_{j\in \Z}{
                                    \abs{c_{j,k}}
                                    \norm{
                                        T_{\sigma_j^k}{
                                            \l(
                                                \Delta_{j+l}^\varphi{f_i}
                                            \r)
                                        }
                                    }_2^2
                                }
                            \r)^\frac{1}{2}
                        }
                    \r)
                }
            }
        }.
    $$
    For $m \in \N_0$ we choose a non-negative Schwarz function $\Phi^m$, such that $\Phi^m(\xi)=1$ for $2^{-7-\ceil{\frac{m}{2}}} \leq \abs{\xi} \leq 2^6$, it vanishes outside of an annulus $2^{-9-\ceil{\frac{m}{2}}} \leq \abs{\xi} \leq 2^8$.
    
    \textbf{Notice}, that $\Delta_{j+l}^\varphi=\Delta_{j+l}^\varphi \Delta_j^{\Phi^m}$ for $l = -5-m, \dots,4$ and supports of $\Phi^m(2^{-j}\cdot)$ overlap at most $m+1$ times up to a multiplicative constant independent of $m$. Also, recall that $\sigma_j^k$ is a Fourier multiplier on $L^2(\R^n)$ with the operator norm less than or equal to 1. The same holds for functions $\varphi(2^{-j}\,\cdot)$. Then
    $$  
        \sum_{j \in \Z}{
            \norm{
                T_{\sigma_j^k}{
                    \l(
                        \Delta_{j+l}^{\varphi}{f_i}
                    \r)
                } 
            }_2^2
        }
        \leq
        \sum_{j \in \Z}{
            \norm{
                \Delta_j^{\Phi^m}{f_i}
            }_2^2
        }
        \lesssim
        (m+1)
        \norm{f_i}_2^2.
    $$

    From this, using Cauchy-Schwarz inequality, Fubini's theorem and Lemma \ref{abs_convergence_lemma}, we have: 
    \begin{align*}
        \norm{T_\sigma^{1,1}(f_1,f_2)}_1
        &\leq
        \sum_{m=0}^\infty{
            \left(\ceil{m/2}+10\right)^2 \,
            \prod_{i=1}^2{
                \l(
                    \sum_{j \in \Z}{
                        \norm{
                            \Delta_j^{\Phi^m}{f_i}
                        }_2^2
                        \sum_{k \in D_m}{
                            \abs{c_{j,k}}
                        }
                    }
                \r)^\frac{1}{2}
            }
        }\\
        &\lesssim
        \sum_{m=0}^\infty{
            (m+1)^2
            \prod_{i=1}^2{
                \l(
                    \sum_{j \in \Z}{
                        \norm{
                            \Delta_j^{\Phi^m}{f_i}
                        }_2^2
                        2^{-\varepsilon m}
                        C
                    }
                \r)^\frac{1}{2}
            }
        } \\
        &\lesssim
        C
        \l(
            \sum_{m=0}^\infty{
                (m+1)^3 \,
                2^{-\varepsilon m}
            }
        \r)
        \norm{f_1}_2
        \norm{f_2}_2.
    \end{align*}
    And because $\varepsilon$ is positive, the sum in the last inequality converges. Hence $T_\sigma^{1,1}$, and consequently $T_\sigma^1$, is bounded from $L^2(\R^n) \times L^2(\R^n)$ to $L^1(\R^n)$.

\subsection{Off-diagonal parts:} 
    
    First we make an observation that there exists a constant $D > 0$, such that for every $m \in \N_0, \, j \in \Z$ and $k \in D_m$, the inequality
    $$
    \abs{
        T_{\sigma_j^k}{
            \l(
                \sum_{r=-\infty}^{j-6-\ceil{\frac{m}{2}}}{
                    \Delta_r^\varphi{f_1}
                }
            \r)
        }
    }
    \leq D \, M(f_1),
    $$
    where $M$ is Hardy-Littlewood maximal operator, holds almost everywhere. 
    
    This is due to the form of the function $\sigma_j^k$  and the fact that the function $\sum_{r=-\infty}^{j-6-\ceil{\frac{m}{2}}}{\varphi(2^{-r} \, \xi)}$ is smooth and has compact support $\abs{\xi} \leq 2^{j-4-\ceil{\frac{m}{2}}}$, where $\ceil{\frac{m}{2}} \approx \log_2(\sqrt{|k|})$. Therefore it is straightforward to check that $\sigma_j^k$ has on this support bounded derivatives independently of $j$ and $k$. 

We have:
    \begin{align*}
        \norm{
            T_\sigma^2(f_1,f_2)
        }_1
        &=
        \norm{
            \sum_{j\in\Z}{
                \sum_{m=0}^\infty{
                    \sum_{k \in D_m}{
                        c_{j,k} \,
                        T_{\sigma_j^k}{
                        \l(
                            \sum_{r=-\infty}^{j-6-\ceil{\frac{m}{2}}}{
                                \Delta_r^\varphi{f_1}
                            }
                        \r)
                        }
                        T_{\sigma_j^k}{
                        \l(
                            \sum_{s=j-4}^{j+4}{
                                \Delta_s^\varphi{f_2}
                            }
                        \r)
                        }   
                    }
                }
            }
        }_1  
    \end{align*}
We fix $m$ and for each $j$ we find a rearrangement $k(j,l)$ $1\leq l \leq 2^{m-1}$ of $D_m$ such that $|c_{j,k(j,l)}|$ is nonincreasing in $l.$As we have $\sum_l |c_{j,k(j,l)}|\leq C2^{-\epsilon m},$ we get $|c_{j,k(j,l)}|\leq C2^{-\epsilon m}/l.$ We also observe that the Fourier support of 
$$T_{\sigma_j^k}	\l(
                            \sum_{r=-\infty}^{j-6-\ceil{\frac{m}{2}}}{
                                \Delta_r^\varphi{f_1}
                            }
                        \r)
                        T_{\sigma_j^k}
                        \l(
                            \sum_{s=j-4}^{j+4}{
                                \Delta_s^\varphi{f_2}
                            }
                        \r)     		
												$$
is contained in an annulus $|x|\in [2^{j-10},2^{j+10}].$ Therefore, we may estimate the $L_1$ norm by the $H_1$ norm expressed by the square function. We get
 \begin{align*}
        &\norm{
            T_\sigma^2(f_1,f_2)
        }_1 \\
        &\leq C
				\sum_{m=0}^\infty\sum_{l=1}^{2^{m-1}}\norm{\l( \sum_{j \in \mathbb Z} \l|c_{j,k(j,l)}T_{\sigma_j^{k(j,l)}}{
                        \l(
                            \sum_{r=-\infty}^{j-6-\ceil{\frac{m}{2}}}
                                \Delta_r^\varphi{f_1}
                        \r)
                        }
                        T_{\sigma_j^{k(j,l)}}
                        \l(
                            \sum_{s=j-4}^{j+4}{
                                \Delta_s^\varphi{f_2}
                            }
                        \r)\r|^2\r)^{1/2}}_1 \\
					&\leq C \sum_{m=0}^\infty\sum_{l=1}^{2^{m-1}}2^{-\epsilon m}/l\norm{ M(f_1)\l( \sum_{j \in \mathbb Z} \l|
                        T_{\sigma_j^{k(j,l)}}
                        \l(
                            \sum_{s=j-4}^{j+4}{
                                \Delta_s^\varphi{f_2}
                            }
                        \r)\r|^2\r)^{1/2}}_1 \\
					&\leq C \norm{f_1}_2 \|f_2\|_2.
	\end{align*}
  
  For the last inequality,we notice that 
    $$
		 \sum_{j \in \mathbb Z} \l|
                        T_{\sigma_j^{k(j,l)}}
                        \l(
                            \sum_{s=j-4}^{j+4}{
                                \Delta_s^\varphi{f_2}
                            }
                        \r)\r|^2 
    $$
		is always sum of almost orthogonal members and the operators $T$ are $L^2$ bounded multipliers. 
\qed

\section{Proof of theorem \ref{Bochner_Riesz}}

Choose two smooth nonnegative radial functions $\varphi_1, \varphi_2 \in C^\infty(\R)$ such that $\varphi_1(\xi) = 0$ for $\abs{x} \leq \frac{1}{4}$, $\varphi_1(\xi) = 1$ for $\abs{x} > \frac{1}{2}$ and $\varphi_1+\varphi_2 \equiv 1$.

For $i \in \{1,2\}$ set $\sigma_i = \sigma^\lambda \, \varphi_i{\l(\abs{\cdot}^2\r)}$. We will split the proof into two parts. First, we will analyze operator $T_{\sigma_1}$ using Lemma \ref{compact_multipliers_lemma}. For $T_{\sigma_2}$ we will check the Hörmander-Mihilin multiplier condition.

\subsection{Boundedness of $T_{\sigma_1}$}

As in Lemma \ref{key_lemma}, we define $\Tilde{\sigma}_1(\pm\abs{\xi}^2) = \sigma_1{\l(\sqrt{2m} \xi\r)}$. We see, that:
$$
    \Tilde{\sigma}_1(x)
    =
    \l( 1 - 2m \abs{x} \r)_+^\lambda \, 
    \varphi_1{ \l(2m \abs{x}\r)}.
$$

Due to the symmetry of $c_k$'s, we consider only $k \in \N$. We have:
\begin{align*}
    c_k 
    =
    \int_{-\frac{1}{2}}^\frac{1}{2}{
        \Tilde{\sigma}_1(x) \,
        e^{-2\pi i k x} \, dx
    }
    &=
    2
    \int_{0}^\frac{1}{2}{
        \l( 1 - 2m x \r)_+^\lambda \, 
        \varphi_1{\l( 2 m x\r)} \,
        \cos{( 2\pi k x)} \, dx
    } \\
    &= \frac{1}{k^{1+\lambda}} \, m^\lambda \,
    \int_0^{\frac{k}{2m}}{
        z^\lambda \,
        \varphi_1{\l(1-\frac{m}{k}z\r)} \,
        \cos{\l( \pi \l(z - \frac{k}{m} \r)\r)} \,
        dz
    }.
\end{align*}

It is not difficult to verify, using integration by parts multiple times, that the integral in the last equality is uniformly bounded with respect to $k$. This shows that $c_k = O(k^{-(1+\lambda)})$ and, consequently, $\{c_k\}_{k \in \Z}$ lies in $\ell^{\frac{2}{m}}(\Z)$, because
$\frac{2}{m} (\lambda+1) > 1$.

Now, using the Lemma \ref{compact_multipliers_lemma}, we get the boundedness of $T_{\sigma_1}$ from $\l(L^2(\R^n)\r)^m$ to $L^\frac{2}{m}(\R^n)$.

\subsection{Boundedness of $T_{\sigma_2}$}

To check the boundedness of the operator $T_{\sigma_2}$, we use the Hörmander-Mihilin multiplier condition (see in \cite{Grafakos_Modern}, Section 7.5.3.).
Fix a Schwartz function $\Phi$ which is supported in an annulus $P$ of the form $\frac{1}{2} \leq \abs{\xi} \leq 2$ and satisfies:
$$
    \sum_{j \in \Z}{\Phi{\l(2^{-j} \xi \r)}} = 1, 
    \quad 
    \forall \xi \in \R^{n} \setminus\{0\}.
$$
We want to show that
$$
    \sup_{k \in \Z}{
        \norm{
            \sigma_2{\l( 2^k \cdot\r)} \Phi
        }_{L^2_{mn}{(P)}}
    }
    <
    \infty,
$$
where $L^2_{mn}(P)$ denotes Sobolev space on $P$, which we will treat in the classical sense, i.e. using weak derivatives.

We see that $\supp{ \sigma_2{\l(2^k\cdot\r)}} \subseteq B\l(0, 2^{-k-\frac{1}{2}} \r)$. This means that $\sigma_2{\l(2^k \cdot\r)} \Phi \equiv 0$ for $k \geq 1$. Hence, we need to focus our attention only on the case of $k$ being a~non-positive integer. Fix~such~$k$.

Using the results of \cite{Ostermann_Radial_Sobolev} we can write, that
\begin{align*}
    \norm{
        \sigma_2{\l(2^k \cdot\r)} \Phi
    }_{L^2_{mn}{(P)}}
    &\leq
    \norm{\Phi}_{L^\infty_{mn}{(P)}}
    \norm{\sigma_2{\l(2^k \cdot\r)}}_{L^2_{mn}(P)} \\
    &\lesssim
    \norm{\Phi}_{L^\infty_{mn}{(P)}}
    \sum_{j=0}^{mn}{
        \norm{
            r^{\frac{mn-2+2j}{4}} \,
            \l(
                \sigma_R\l(2^k r\r)
            \r)^{(j)}
        }_{L^2(2^{-2},4)}
    },
\end{align*}
where $\sigma_R(r) = (1-r)_+^\lambda \varphi(r)$. It holds that $\sigma_R$ is a compactly supported smooth function and $\sigma_R{\l(\abs{\xi}^2\r)} = \sigma_2(\xi)$.

For each $j \in \{0, \dots, mn\}$, we have
\begin{align*}
    \norm{
            r^{\frac{mn-2+2j}{4}} \,
            \l(
                \sigma_R\l(2^k r\r)
            \r)^{(j)}
    }_{L^2(2^{-2},4)}
    &=
    2^{jk} 
    \norm{
            r^{\frac{mn-2+2j}{4}} \,
            \sigma_R^{(j)}\l(2^k r\r)
    }_{L^2(2^{-2},4)} \\
    &= 
    2^{k\l(j-\frac{d+2j}{4} \r)}
    \norm{
            r^{\frac{mn-2+2j}{4}} \,
            \sigma_R^{(j)}\l(r\r)
    }_{L^2(2^{k-2},2^{k+2})} \\
    &\lesssim
    2^{jk}
    \norm{\sigma_R}_{L^\infty_{mn}(\R)} \\
    &\leq 
    \norm{\sigma_R}_{L^\infty_{mn}(\R)}.
\end{align*}
It follows, that
$$
    \norm{
        \sigma_2{\l(2^k \cdot\r)} \Phi
    }_{L^2_{mn}{(P)}}
    \lesssim
    (mn+1)
    \norm{
        \Phi
    }_{L^\infty_{mn}(P)}
    \norm{
        \sigma_R
    }_{L^\infty_{mn}(\R)} 
    < \infty.
$$
So $T_{\sigma_2}$ is bounded from $\l(L^2 \l(\R^n\r)\r)^m$ to $L^{\frac{2}{m}}(\R^n)$.
\qed

\section{Proof of theorem \ref{Modified_Bochner_Riesz}}

The proof is analogous to the first part of the proof of theorem \ref{Bochner_Riesz} for the case $m=2$. Again, we take $\Tilde{\sigma}^{0,\gamma}(\pm\abs{\xi}^2) = \sigma^{0, \gamma}(2\xi)$. Then for $k \in \N$, we have
$$
    c_k 
    =
    \int_{-\frac{1}{2}}^\frac{1}{2}{
        \Tilde{\sigma}^{0, \,\gamma}(x) \,
        e^{-2\pi k i x} \, dx
    }
    =
    \frac{1}{2}
    \int_0^1{
        \l(1-\log{\l(1-x\r)}\r)^{-\gamma}
        \cos{
            \l(
                \frac{\pi}{2}k x
            \r)
        }
        \, dx
    }.
$$

We have $c_k = O{\l( \abs{k}^{-1} \log{\l(\abs{k}\r)}^{-\gamma}\r)}$. To verify it, we will split the~integral into 2 parts. One on the interval $\left[1-\frac{1}{k}, 1\right]$ and the second on $\left[0, 1-\frac{1}{k}\right]$.

The first part of the integral is straightforward:
$$
    \int_{1-\frac{1}{k}}^1{
        \abs{
            \l(1-\log{\l(1-x\r)}\r)^{-\gamma}
            \cos{
                \l(
                    \frac{\pi}{2}k x
                \r)
            }
        }
        \, dx
    }
    \leq
    \frac{1}{k}
    \l(1+\log{\l(k\r)}\r)^{-\gamma}
    \approx
    \frac{1}{k \log{\l(k\r)}^\gamma}.
$$
The second part requires integration by parts 2 times. For sufficiently large $k$'s we obtain:
\begin{align*}
    \abs{
        \int_0^{1-\frac{1}{k}}{
            \l(1-\log{\l(1-x\r)}\r)^{-\gamma}
                \cos{
                    \l(
                        \frac{\pi}{2}k x
                    \r)
                }
            \, dx
        } 
    }
    &\lesssim
    \frac{1}{k \log{\l(k\r)}^\gamma}
    +
    \frac{1}{k \log{\l(k\r)}^{\gamma+1}}
    +
    \frac{1}{k^2} \\
    &\lesssim
    \frac{1}{k \log{\l(k\r)}^\gamma},
\end{align*}
Therefore, truly $c_k = O{\l( \abs{k}^{-1} \log{\l(\abs{k}\r)}^{-\gamma}\r)}$.

Then $\{c_k\}_{k \in \Z}$ lies in $\ell^1(\Z)$, because $\gamma > 1$. So, again by Lemma \ref{compact_multipliers_lemma}, we get the boundedness of $T_{\sigma^{0, \gamma}}$ from $L^2(\R^n) \times L^2(\R^n)$ to $L^1(\R^n)$.
\qed

\bibliographystyle{plain}

\renewcommand{\bibname}{Literature}

\bibliography{literatura}

\end{document}